\documentclass[11pt]{amsart}
\usepackage{graphicx}
\usepackage{amssymb}

\newcommand{\R}{\mathbb{R}}
\newcommand{\N}{\mathbb{N}}
\newcommand{\half}{\frac{1}{2}}

\let\Re=\undefined\DeclareMathOperator{\Re}{Re}
\let\Im=\undefined\DeclareMathOperator{\Im}{Im}

\def\C{{\mathbb C}}
\def\R{{{\mathbb R}}}
\def\Z{{{\mathbb Z}}}

\newcommand{\eps}{\varepsilon}

\newcommand{\ir}{I\times\R}

\newcommand{\zit}{Z_I(t)}

\newcommand{\itxi}{\int_{t_0}^t\int_{\sum_{i=1}^{2k+2}\xi_i=0}}
\newcommand{\symb}
{1-\frac{m(\xi_2+\xi_3+\cdots+\xi_{2k+2})}{m(\xi_2)m(\xi_3)\cdots m(\xi_{2k+2})}}

\newcommand{\ma}{L_{t,x}^8}

\newcommand{\hs}{H_x^s}
\newcommand{\dhs}{\dot H_x^s}
\newcommand{\ho}{\dot H_x^1}
\newcommand{\ulam}{u^{\lambda}}
\newcommand{\lxt}{L_x^2}

\newcommand{\tk}{\frac 12-\frac 1k}

\newcommand{\ztr}{[0,T]\times\R}

\newtheorem{theorem}{Theorem}
\theoremstyle{definition}
\newtheorem{definition}{Definition}
\theoremstyle{remark}

\theoremstyle{proposition}
\newtheorem{proposition}{Proposition}
\theoremstyle{lemma}
\newtheorem{lemma}{Lemma}
\theoremstyle{corollary}

\numberwithin{equation}{section}
\numberwithin{lemma}{section}
\numberwithin{remark}{section}
\numberwithin{theorem}{section}
\numberwithin{corollary}{section}
\numberwithin{proposition}{section}
\numberwithin{definition}{section}

\title[Scattering for generalized NLS on $\R$]{Global existence and scattering
for rough solutions to generalized nonlinear Schr\"odinger equations on $\R$}
\author{J.~Colliander}
\address{University of Toronto and M.S.R.I.}
\author{J.~Holmer}
\address{University of California, Berkeley and M.S.R.I.}
\author{M.~Visan}
\address{Institute for Advanced Study and M.S.R.I.}
\author{X.~Zhang}
\address{Academy of Mathematics and System Sciences, Chinese Academy of Sciences and M.S.R.I.}


\begin{document}

\begin{abstract}
We consider the Cauchy problem for a family of semilinear defocusing Schr\"odinger
equations with monomial nonlinearities in one space dimension.  We establish global
well-posedness and scattering.  Our analysis is based on a four-particle interaction Morawetz
estimate giving \emph{a priori} $L^8_{t,x}$ spacetime control on solutions.
\end{abstract}

\maketitle

\section{Introduction}
We consider the initial value problem for the one-dimensional defocusing nonlinear Schr\"odinger (NLS) equation,
\begin{equation}
\label{nls}
\left\{
\begin{matrix}
i u_t + \Delta u    = |u|^{2k} u  \\
 u(0,x) = u_0(x),
\end{matrix}
\right.
\end{equation}
where $k \in \N$ with $k \geq 3$ and $u$ is a complex-valued function on spacetime $\R_t\times\R_x$.
This problem is known to be locally wellposed for initial data in $H^s (\R)$ for $s \geq s_c : = \half - \frac{1}{k}$; see \cite{cw:local, cazbook}.
The scaling invariant Sobolev index $s_c$ is distinguished in the theory by the invariance of the $\dot H_x^{s_c}$ norm
under the scaling symmetry of solutions to \eqref{nls}:  If $u$ solves \eqref{nls} then
\begin{equation}
\label{scaling}
u^\lambda (t ,x ) := \lambda^{-\frac{1}{k}} u(\lambda^{-2} t , \lambda^{-1} x)
\end{equation}
also solves \eqref{nls}.

The following quantities, if finite for the initial data, are time invariant:
\begin{equation*}
\begin{split}
\text{Mass}&:=M[u(t)] := \|u(t)\|_{\lxt}^2,\\
\text{Energy}&:=E[u(t)] := \half \|\nabla u(t)\|_{\lxt}^2 + \frac{1}{2k+2} \|u(t)\|_{L_x^{2k+2}}^{2k+2}.
\end{split}
\end{equation*}
The local-in-time theory in the presence of these conserved quantities iterates to prove global-in-time well-posedness for \eqref{nls} for initial
data in $H_x^1$. Furthermore, in this case it is known that these global-in-time solutions are bounded in the associated scaling-invariant
diagonal Strichartz space $L^{3k}_{t,x}$ and scatter; see \cite{Nakanishi}.  It is conjectured  that global well-posedness and scattering
also hold for solutions to \eqref{nls} with initial data in $\dot H^{s_c}(\R)$.

This work makes partial progress toward this conjecture by establishing
these properties for solutions to \eqref{nls} with initial data in $H^{\frac{8}{9} } (\R)$.  In fact, for all values $k$ considered
we establish global well-posedness and scattering for \eqref{nls} with initial data in $H^s (\R)$ for $s> s_k$, where
$s_k := \frac{8k - 16}{9k - 14} < \frac{8}{9}$.

\begin{theorem}\label{main}
For each $k \in \{ 3, 4, \dots \}$ there is a regularity threshold $s_k = \frac{8k - 16}{9k -14}$ such that the initial
value problem \eqref{nls} is globally wellposed and scatters for
initial data $u_0 \in H^{s} (\R)$, provided $s > s_k$. In particular,
there exist $u_{\pm} \in H^s (\R)$ such that
\begin{equation*}
  \| u(t) - e^{it \Delta} u_{\pm} \|_{H^s (\R)} \to 0 \quad \text{as }
  t \to \pm \infty.
\end{equation*}
\end{theorem}

Our approach to proving this result is based on the proof of a similar statement for the defocusing cubic nonlinear
Schr\"odinger equation on $\R^3$ in \cite{CKSTT:CPAM}.  The analysis in \cite{CKSTT:CPAM} is based on an \emph{a priori}
two-particle interaction Morawetz estimate.  We derive a four-particle interaction Morawetz inequality
which provides $L^8_{t,x}$ spacetime control on solutions to \eqref{nls}.  Our analysis relies on this improved \emph{a priori} control.

As a consequence of the four-particle interaction Morawetz inequality, we are in fact able to offer a new proof of scattering
for a class of one-dimensional defocusing nonlinear Schr\"odinger equations with initial data in $H^1(\R)$; see \cite{Nakanishi}
for the original proof.

\begin{theorem}[Scattering in $H^1(\R)$]\label{scat}
Let $u_0\in H^1(\R)$.  Then, there exists a unique global solution $u$ to the initial value problem
\begin{equation}\label{nls scat}
\begin{cases}
i u_t +\Delta u = |u|^{2p}u, \quad p>0,\\
u(0,x) = u_0(x).
\end{cases}
\end{equation}
Moreover, if $p>2$ there exist $u_\pm\in H^1(\R)$ such that
$$
\|u(t)-e^{it\Delta}u_\pm\|_{H^1(\R)}\to 0 \quad \text{as } t\to \pm \infty.
$$
\end{theorem}

We briefly explain our strategy for proving Theorem~\ref{main} and Theorem~\ref{scat}.

The interaction Morawetz inequality we derive in Section~\ref{MorawetzSection} provides \emph{a priori} $L^8_{t,x}$ spacetime control
on solutions to \eqref{nls scat} (and hence on solutions to \eqref{nls}), provided that $\| u(t) \|_{H_x^{1/2}}$ stays bounded.
In particular, if the initial data $u_0\in H_x^1$, we immediately obtain that the unique global $H_x^1$ solution enjoys
the global $\ma$ estimate.  In Section~\ref{S:scat}, for $p>2$ we upgrade this estimate to stronger Strichartz norm control from which
scattering in $H_x^1$ follows, thus establishing Theorem~\ref{scat}.  A similar argument in higher dimensions, $n\ge 3$,
relying on the two-particle Morawetz inequality, can be found in \cite{tvz}.

If we are in the $\hs$ setting (rather than the $H^1_x$ setting) with $s$ being defined in Theorem~\ref{main},
we know the problem is $\hs$ subcritical and, as a consequence, the length of the local well-posedness time interval
of the unique $\hs$ solution depends only on the $\hs$ norm of the initial data.  Thus, in order to prove global well-posedness
we only need to control the $\hs$ norm of the solution.  This is not immediate as the $H^s_x$ norm is not conserved.
In order to derive the desired control over the $\hs$ norm of the solution, we will use the `$I$-method'.

The idea behind the `$I$-method' (\cite{CKSTT:MRL, CKSTT:CPAM}) is to smooth out the initial data in
order to get access to the good local and global theory available at $H_x^1$ regularity.  To this end, one introduces the Fourier
multiplier $I$ which is the identity on low frequencies and behaves like a fractional integral operator of order $1-s$ on high
frequencies. Thus, the operator $I$ maps $\hs$ to $H_x^1$ and the $\hs$ norm of $u$ can be controlled by the $H_x^1$ norm of
the modified solution $Iu$.  However, $Iu$ is not a solution to \eqref{nls} and hence one cannot use the conservation of energy
to derive a bound on the $H^1_x$ norm of $Iu$.  In fact, we expect an increment in the energy of $Iu$.  This increment is proved to be
under control provided the Morawetz norm is finite; see Section~\ref{S:ac}.  But in order for the Morawetz norm to be finite we need
to control the $H^{1/2}_x$ norm of the solution.  This sets us up for a bootstrap argument which will be carried out in Section~\ref{S:boot}.

\textbf{Acknowledgments:}
J.~Colliander is partially supported by N.S.E.R.C. Grant RGPIN 250233-03.  J.~Holmer is partially supported by an N.S.F. postdoctoral fellowship.
M.~Visan is supported by the N.S.F. grant DMS 0111298.  X.~Zhang is supported by the NSF grant No. 10601060 (China).
We thank Terry Tao for useful discussions related to this work. We gratefully acknowledge support from the Mathematical Sciences
Research Institute where this work was completed.

\section{Preliminaries}
In this section, we introduce notations and some basic estimates we will invoke throughout this paper.

We will often use the notation $X \lesssim Y$ whenever there exists some constant $C>0$ so that $X \leq CY$.  Similarly, we will use $X \sim Y$ if $X
\lesssim Y \lesssim X$.  We will use $X\ll Y$ if $X\leq cY$ for some very small constant $c>0$.
We will sometimes denote partial derivatives with subscripts ($a_j(x):=\partial_j a(x):=\partial_{x_j} a(x)$)
and use the convention that repeated indices are implicitly summed.

We use $L_x^r(\R)$ to denote the Banach space of functions $f:\R\to \C$ whose norm
$$
\|f\|_r:=\Bigl(\int_{\R} |f(x)|^r dx\Bigr)^{1/r}
$$
is finite, with the usual modifications when $r=\infty$.

We use $L_t^qL_x^r$ to denote the spacetime norm
$$
\|u\|_{q,r}:=\|u\|_{L_t^qL_x^r(\R\times\R)}
:=\Bigl(\int_{\R}\Bigl(\int_{\R} |u(t,x)|^r dx \Bigr)^{q/r}dt\Bigr)^{1/q},
$$
with the usual modifications when either $q$ or $r$ are infinity, or when the domain $\R \times \R$ is replaced by some smaller
spacetime region.  When $q=r$ we abbreviate $L_t^qL_x^r$ by $L^q_{t,x}$.

We define the Fourier transform on $\R$ to be
$$
\hat f(\xi) := \int_{\R} e^{-2 \pi i x \cdot \xi} f(x) dx.
$$

We will make use of the fractional differentiation operators $|\nabla|^s$ defined by
$$
\widehat{|\nabla|^sf}(\xi) := |\xi|^s \hat f (\xi).
$$
These define the homogeneous Sobolev norms
$$
\|f\|_{\dot H^{s}_x} := \| |\nabla|^s f \|_{L^2_x}
$$
and more general Sobolev norms
$$
\|f\|_{H_x^{s,p}}:=\|\langle\nabla\rangle^s f\|_p,
$$
where, $\langle\nabla\rangle=(1+|\nabla|^2)^{\frac 12}$.

Let $e^{it\Delta}$ be the free Schr\"odinger propagator.  In physical space this is given by the formula
$$
e^{it\Delta}f(x) = \frac{1}{(4 \pi i t)^{1/2}} \int_{\R} e^{i|x-y|^2/4t} f(y) dy
$$
for $t\neq 0$ (using a suitable branch cut to define $(4\pi it)^{1/2}$), while in frequency space one can write this as
\begin{equation}\label{fourier rep}
\widehat{e^{it\Delta}f}(\xi) = e^{-4 \pi^2 i t |\xi|^2}\hat f(\xi).
\end{equation}
In particular, the propagator obeys the \emph{dispersive inequality}
\begin{equation}\label{dispersive ineq}
\|e^{it\Delta}f\|_{L^\infty_x} \lesssim |t|^{-\frac{1}{2}}\|f\|_{L^1_x}
\end{equation}
for all times $t\neq 0$.

We also recall \emph{Duhamel's formula}
\begin{align}\label{duhamel}
u(t) = e^{i(t-t_0)\Delta}u(t_0) - i \int_{t_0}^t e^{i(t-s)\Delta}(iu_t + \Delta u)(s) ds.
\end{align}

\begin{definition}
A pair of exponents $(q,r)$ is called \emph{Schr\"odinger-admissible} if
$$
\frac{2}{q} +\frac{1}{r} = \frac{1}{2}, \quad  2 \leq r \leq \infty.
$$
\end{definition}
For a spacetime slab $I\times\R$, we define the Strichartz norm
$$
\|f\|_{S^0(I)}:=\sup_{(q,r)\text{ admissible}}\|f\|_{L_t^qL_x^r(I\times\R)}.
$$
Then, we have the following Strichartz estimates (for a proof see \cite{gv:strichartz, tao:keel,strichartz}):

\begin{lemma}\label{lemma linear strichartz}
Let $I$ be a compact time interval, $t_0\in I$, $s\ge 0$, and let $u$ be a solution to the forced Schr\"odinger equation
\begin{equation*}
i u_t + \Delta u =\sum_{i=1}^m F_i
\end{equation*}
for some functions $F_1, \dots, F_m$.  Then,
\begin{equation}
\||\nabla|^s u\|_{ S^0(I)} \lesssim \|u(t_0)\|_{\dhs} + \sum_{i=1}^m \||\nabla| ^s F_i\|_{L_t^{q_i'}L_x^{r_i'}(I\times\R)}
\end{equation}
for any admissible pairs $(q_i,r_i)$, $1\leq i\leq m$.  Here, $p'$ denotes the conjugate exponent to $p$, that is, $\tfrac 1p + \tfrac 1{p'}=1$.
\end{lemma}

We will also need some Littlewood-Paley theory.  Specifically, let $\varphi(\xi)$ be a smooth bump supported in $|\xi| \leq 2$
and equalling one on $|\xi| \leq 1$.  For each dyadic number $N \in 2^\Z$ we define the Littlewood-Paley operators
\begin{align*}
\widehat{P_{\leq N}f}(\xi) &:=  \varphi(\xi/N)\hat f (\xi),\\
\widehat{P_{> N}f}(\xi) &:=  [1-\varphi(\xi/N)]\hat f (\xi),\\
\widehat{P_N f}(\xi) &:=  [\varphi(\xi/N) - \varphi (2 \xi /N)] \hat f (\xi).
\end{align*}
Similarly, we can define $P_{<N}$, $P_{\geq N}$, and $P_{M < \cdot \leq N} := P_{\leq N} - P_{\leq M}$, whenever $M$ and $N$ are dyadic
numbers.  We will frequently write $f_{\leq N}$ for $P_{\leq N} f$ and similarly for the other operators. We recall the following
standard Bernstein and Sobolev type inequalities:

\begin{lemma}\label{bernstein}
For any $1\le p\le q\le\infty$ and $s>0$, we have
\begin{align*}
\|P_{\geq N} f\|_{L^p_x} &\lesssim N^{-s} \| |\nabla|^s P_{\geq N} f \|_{L^p_x}\\
\| |\nabla|^s  P_{\leq N} f\|_{L^p_x} &\lesssim N^{s} \| P_{\leq N} f\|_{L^p_x}\\
\| |\nabla|^{\pm s} P_N f\|_{L^p_x} &\sim N^{\pm s} \| P_N f \|_{L^p_x}\\
\|P_{\leq N} f\|_{L^q_x} &\lesssim N^{\frac{1}{p}-\frac{1}{q}} \|P_{\leq N} f\|_{L^p_x}\\
\|P_N f\|_{L^q_x} &\lesssim N^{\frac{1}{p}-\frac{1}{q}} \| P_N f\|_{L^p_x}.
\end{align*}
\end{lemma}

\vspace{0.4cm}

For $N>1$, we define the Fourier multiplier $I:=I_N$ (cf. \cite{CKSTT:MRL})
$$
\widehat{I_N u}(\xi):=m_N(\xi)\hat u(\xi),
$$
where $m_N$ is a smooth radial decreasing function such that
$$
m_N(\xi)=\left\{
\begin{array}{cc}
1, & \text{if} \quad  |\xi|\le N\\
\bigl(\frac{|\xi|}{N}\bigr)^{s-1}, & \text{if} \quad |\xi|\ge
2N.
\end{array}
\right.
$$
Thus, $I$ is the identity operator on frequencies $|\xi|\le N$ and behaves like a fractional integral operator of order $1-s$ on higher
frequencies. In particular, $I$ maps $H^s_x$ to $H_x^1$.  We collect the basic properties of $I$ into the following

\begin{lemma}\label{basic property}
Let $1<p<\infty$ and $0\leq\sigma\le s<1$.  Then,
\begin{align}
\|If\|_p&\lesssim \|f\|_p \label{i1}\\
\||\nabla|^\sigma P_{> N}f\|_p&\lesssim N^{\sigma-1}\|\nabla I f\|_p \label{i2}\\
\|f\|_{H^s_x}\lesssim \|If\|_{H^1_x}&\lesssim
N^{1-s}\|f\|_{H^s_x}.\label{i3}
\end{align}
\end{lemma}
\begin{proof}
The estimate \eqref{i1} is a direct consequence of the multiplier
theorem.

To prove \eqref{i2}, we write
$$
\||\nabla|^\sigma P_{> N} f\|_p=\|P_{> N}|\nabla |^\sigma(\nabla
I)^{-1}\nabla I f\|_p.
$$
The claim follows again from the multiplier theorem.

Now we turn to \eqref{i3}.  By the definition of the operator $I$
and \eqref{i2},
\begin{align*}
\|f\|_{H^s_x}&\lesssim \|P_{\le N} f\|_{H^s_x}+\|P_{>N}f\|_2+\||\nabla|^s P_{>N} f\|_2\\
&\lesssim \|P_{\le N} I f\|_{H_x^1}+N^{-1}\|\nabla I f\|_2+N^{s-1}\|\nabla I f\|_2\\
&\lesssim \|If\|_{H^1_x}.
\end{align*}
On the other hand, since the operator $I$ commutes with $\langle
\nabla\rangle^s$,
\begin{align*}
\|If\|_{H_x^1} =\|\langle\nabla\rangle^{1-s}I\langle\nabla\rangle^s
f\|_2 \lesssim N^{1-s}\|\langle\nabla\rangle^sf\|_2 \lesssim
N^{1-s}\|f\|_{H^s_x},
\end{align*}
which proves the last inequality in \eqref{i3}.  Note that a similar
argument also yields
\begin{align}\label{i4}
\|If\|_{\dot H^1_x}&\lesssim N^{1-s}\|f\|_{\dot H^s_x}.
\end{align}
\end{proof}

\section{An interaction Morawetz inequality}\label{MorawetzSection}
In this section we develop an \emph{a priori} four-particle interaction Morawetz inequality for solutions to one-dimensional defocusing
nonlinear Schr\"odinger equations.  This \emph{a priori} control will be fundamental to our analysis.

The name Morawetz inequality derives from her work on monotonicity formulae for the wave equation.  The Schr\"odinger version is
due to Lin and Strauss, \cite{linstrauss}.  The idea of a two-particle interaction Morawetz inequality was first introduced in \cite{CKSTT:CPAM}.
This two-particle style of estimate has proved invaluable in the study of NLS in dimensions three and higher.  Unfortunately, there is no
direct analogue of this estimate in dimensions one and two; nevertheless, several alternatives have been proposed, \cite{Nakanishi, fanggrillakis}.
Here we derive a Morawetz inequality based on four-particle interactions.  This approach was suggested to us by Terry Tao, based on a private
conversation with Andrew Hassel.

\begin{proposition}[Interaction Morawetz estimate]\label{fpim}
Let $u$ be an $H^{1/2}$ solution to \eqref{nls scat} on the spacetime slab $I\times\R$.  Then,
\begin{align}\label{fpim est}
\int_I\int_{\R} |u(t,x)|^8\,dx\,dt \lesssim \|u\|_{L_t^\infty \dot H^{1/2}_x(I\times\R)}^2 \|u_0\|_2^6.
\end{align}
\end{proposition}

The calculations that follow are difficult to justify without additional regularity and decay assumptions on the solution. This
obstacle can be dealt with in the standard manner: mollify the initial data and the nonlinearity to make the interim
calculations valid and observe that the mollifications can be removed at the end. For expository reasons, we skip
the details and keep all computations on a formal level.

In order to prove Proposition~\ref{fpim} we first review general facts about the one-particle Morawetz action.
Let $\phi:\R_t\times\R_y^4\to \C$ be a solution to the Schr\"odinger equation
$$
i\phi_t + \Delta \phi = \mathcal N.
$$
Let $a:\R^4_y\to \R$ be a convex weight function and define the Morawetz action to be the weighted momentum
$$
M_a(t):=2\Im \int_{\R^4}\overline{\phi(t,y)}\nabla a(y)\cdot \nabla \phi(t,y)\,dy.
$$
A direct calculation establishes that in the $(y_1,\dots, y_4)$ coordinate system we have
\begin{align*}
\partial_t M_a(t)
&=2\int_{\R^4}(-\Delta\Delta a(y))|\phi(t,y)|^2\, dy + 4\int_{\R^4}a_{jk}(y)\Re(\overline{\phi_j}\phi_k)(t,y)\, dy\\
&\quad +2\int_{\R^4}\nabla a(y)\cdot \{\mathcal N,\phi\}(t,y)\,dy,
\end{align*}
where the momentum bracket is defined by
$$
\{f,g\}:=\Re (f \nabla \bar g - g \nabla \bar f).
$$
As the weight $a$ is convex, the matrix $\{a_{jk}\}_{1\leq j,k\leq 4}$ is positive semi-definite and hence
$$
\int_{\R^4}a_{jk}(y)\Re(\overline{\phi_j}\phi_k)(t,y)\, dy\geq 0.
$$
Thus,
\begin{align}
\partial_t M_a(t)
&\geq 2\int_{\R^4}(-\Delta\Delta a(y))|\phi(t,y)|^2\, dy \label{ma}\\
&\quad + 2\int_{\R^4}\nabla a(y)\cdot\{\mathcal N,\phi\}(t,y)\,dy.\notag
\end{align}

Now we are ready to prove Proposition~\ref{fpim}.  Let $u$ be a solution to \eqref{nls scat} and for each $1\leq j\leq 4$ let $u_j(t,x_j):=u(t,x_j)$.
Define
$$
w(t,x)=w(t,x_1,x_2,x_3,x_4):=\prod_{j=1}^4 u_j(t,x_j);
$$
note that $w$ satisfies the equation
$$
iw_t+\Delta_x w=\bigl(\sum_{j=1}^4|u_j|^{2p}\bigr)w.
$$
Next, we perform the orthonormal change of variables
$$
z=Ax \mbox{ with } A=\frac 12
\begin{bmatrix}
1 & 1 & 1 & 1 \\
1 & 1 & -1 & -1 \\
1 & -1 & 1 & -1 \\
-1 & 1 & 1 & -1 \\
\end{bmatrix}.
$$
Then $\Delta_x=\Delta_z$ and hence, for $\omega(t,z):=w(t,x(z))$, we have
$$
i\omega_t+\Delta_z \omega=\bigl(\sum_{j=1}^4|u_j|^{2p}\bigr)\omega.
$$
Applying \eqref{ma} to $\omega$ in the $(z_1,\dots, z_4)$ coordinate system with the convex weight $a(z):=(z_2^2+z_3^2+z_4^2)^{1/2}$, we get
\begin{align}
\partial_t M_a(t)
&\geq 2\int_{\R^4}(-\Delta_z\Delta_z a(z))|\omega(t,z)|^2\, dz \notag\\
&\quad + 2\int_{\R^4}\nabla_z a(z)\cdot \bigl\{\bigl(\sum_{j=1}^4|u_j|^{2p}\bigr)\omega,\omega\bigr\}(t,z)\,dz,\label{ma omega}
\end{align}
where
$$
M_a(t):=2\Im \int_{\R^4} \overline{\omega(t,z)} \nabla_z a(z) \cdot \nabla_z \omega(t,z)\, dz.
$$

A quick computation shows that
$$
-\Delta_z\Delta_z a(z)=4\pi \delta(z_2,z_3,z_4)
$$
and hence, by a change of variables,
\begin{align*}
2\int_{\R^4}(-\Delta_z\Delta_z a(z))|\omega(t,z)|^2\, dz_1
&= 8 \pi \int_{\R} |\omega(t,z_1,0,0,0)|^2 \,dz_1\\
&= 16 \pi \int_{\R} |w(t,z_1,z_1,z_1,z_1)|^2 \,dz_1\\
&= 16\pi \int_{\R} |u(t,z_1)|^8\, dz_1.
\end{align*}

To estimate the second term on the right-hand side of \eqref{ma omega}, we note that orthonormal changes of variables leave
inner products invariant and hence,
\begin{align*}
\int_{\R^4}\nabla_z a(z)\cdot \bigl\{\bigl(\sum_{j=1}^4|u_j|^{2p}\bigr)& \omega,\omega\bigr\}(t,z)\,dz\\
&=\int_{\R^4}\nabla_x a(x)\cdot \bigl\{\bigl(\sum_{j=1}^4|u_j|^{2p}\bigr)w,w\bigr\}(t,x)\,dx.
\end{align*}
A simple computation then shows that in the $(x_1, \dots, x_4)$ coordinate system we have
\begin{align*}
\bigl\{\bigl(\sum_{j=1}^4|u_j|^{2p}\bigr)w,w\bigr\}^i
&= \bigl(\sum_{j=1}^4|u_j|^{2p}\bigr)w\partial_{x_i} \bar w - w\partial_{x_i}\bigl[\bigl(\sum_{j=1}^4|u_j|^{2p}\bigr)\bar w\bigr]\\
&= -|w|^2 \partial_{x_i} \bigl(\sum_{j=1}^4|u_j|^{2p}\bigr)\\
&= -\frac{p}{p+1} \partial_{x_i} \bigl(|w|^2 |u_i|^{2p} \bigr).
\end{align*}
Integrating by parts, we obtain
\begin{align*}
\int_{\R^4}\nabla_x a(x)\cdot &\bigl\{\bigl(\sum_{j=1}^4|u_j|^{2p}\bigr)w,w\bigr\}(t,x)\,dx\\
&=\frac{p}{p+1} \int_{\R^4}\sum_{i=1}^4 a_{ii}(x)\bigl(|w|^2 |u_i|^{2p} \bigr)(t,x)\, dx \geq 0,
\end{align*}
as $a$ is a convex function.

Putting everything together we get
$$
\partial_t M_a(t) \geq 8\pi \int_{\R} |u(t,x)|^8\, dx
$$
and hence, by the Fundamental Theorem of Calculus,
\begin{align*}
\int_I \int_{\R} |u(t,x)|^8\, dx\, dt \lesssim \sup_{t\in I} |M_a(t)|.
\end{align*}

In order to estimate the right-hand side in the inequality above, we first note that
\begin{align}\label{hardy}
\Bigl|\int_{\R^n}f(x) \frac{x}{|x|}\cdot \nabla f(x)\, dx\Bigr|\lesssim \|f\|_{\dot H^{1/2}(\R^n)}^2,
\end{align}
for any function $f:\R^n\to \C$ with $n\geq 3$.  Indeed, by Cauchy-Schwarz,
\begin{align*}
\Bigl|\int_{\R^n}f(x) \frac{x}{|x|}\cdot \nabla f(x)\, dx\Bigr|
\lesssim \|f\|_{\dot H^{1/2}(\R^n)} \Bigl\| \frac{x}{|x|} f \Bigr\|_{\dot H^{1/2}(\R^n)},
\end{align*}
and \eqref{hardy} follows if we establish that the operator $T(f)(x):=\frac{x}{|x|} f(x)$ is bounded on $\dot H^{1/2}(\R^n)$.
Using Hardy's inequality
$$
\Bigl\|\frac{f}{|x|}\Bigr\|_2\lesssim \|\nabla f\|_2,
$$
it is easy to see that $T$ is bounded on $L^2(\R^n)$ and on $\dot H^1(\R^n)$.  By interpolation, this yields the claim.

Applying \eqref{hardy} (in the variables $(z_2,z_3,z_4)$), Plancherel, and a change of variables, we estimate
\begin{align*}
|M_a(t)|
&\lesssim \int_\R \|\omega(t,z_1, \cdot)\|_{\dot H^{1/2}(\R^3)}^2\, dz_1\\
&= \int_\R \int_{\R^3} |\xi_2^2+\xi_3^2 +\xi_4^2|^{1/2}|\tilde \omega(t,z_1,\xi_2,\xi_3,\xi_4)|^2\, d\xi_2\, d\xi_3\, d\xi_4\, dz_1\\
&= \int_{\R^4} |\xi_2^2+\xi_3^2 +\xi_4^2|^{1/2}|\hat \omega(t,\xi)|^2\, d\xi\\
&\leq \int_{\R^4} |\xi||\hat \omega(t,\xi)|^2\, d\xi\\
&= \int_{\R^4} |\eta||\hat w(t,\eta)|^2\, d\eta\\
&= \int_{\R^4} |\eta||\hat u_1(t,\eta_1)|^2 |\hat u_2(t,\eta_2)|^2 |\hat u_3(t,\eta_3)|^2 |\hat u_4(t,\eta_4)|^2\, d\eta\\
&\leq \int_{\R^4} \bigl(|\eta_1|+|\eta_2|+|\eta_3|+|\eta_4|\bigr)\prod_{j=1}^4 |\hat u_j(t,\eta_j)|^2\, d\eta\\
&\leq 4 \|u(t)\|_{\dot H^{1/2}}^2 \|u(t)\|_2^6.
\end{align*}
In the computations above,  we used $\tilde \omega$ to denote the partial Fourier transform with respect to the variables $(z_2,z_3,z_4)$ and
$\hat \omega$ to denote the full Fourier transform.  The change of variables performed was $\xi:=A\eta$.

Thus, by the conservation of mass,
\begin{align*}
\int_I \int_{\R} |u(t,x)|^8\, dx\, dt
\lesssim \sup_{t\in I} \|u(t)\|_{\dot H_x^{1/2}}^2 \|u(t)\|_2^6
\lesssim \|u\|_{L_t^\infty \dot H_x^{1/2}(I\times\R)}^2 \|u_0\|_2^6.
\end{align*}
This concludes the proof of Proposition~\ref{fpim}.

\section{Proof of Theorem~\ref{scat}}\label{S:scat}
In this section we prove Theorem~\ref{scat}.  Global well-posedness for \eqref{nls scat} is a consequence of the fact that the equation is
subcritical with respect to energy.  The result and the proof are by now standard and we will not revisit them here; see \cite{cazbook, kato}.

Scattering in the case $p>2$ was first proved by Nakanishi, \cite{Nakanishi}.  In this section we present a new proof
relying on the four-particle interaction Morawetz inequality we developed in the previous section.

Indeed, by Proposition~\ref{fpim} and the conservation of mass and energy, the unique global solution to \eqref{nls scat}
with initial data in $H^1(\R)$ satisfies
\begin{align}\label{mb}
\|u\|_{L_{t,x}^8(\R\times\R)}\lesssim \|u_0\|_{H^1(\R)}.
\end{align}
In order to prove scattering, we first upgrade \eqref{mb} to Strichartz control.
Let $\delta>0$ be a small constant to be chosen momentarily
and divide $\R$ into $L=L(\|u_0\|_{H^1(\R)})$ subintervals $I_j=[t_j, t_{j+1}]$ such that
\begin{align}\label{small it}
\|u\|_{L_{t,x}^8(I_j\times\R)}\sim \delta.
\end{align}
As $p>2$, there exists $\eps>0$ such that $p>2+\frac 18 \eps$.
By Lemma~\ref{lemma linear strichartz}, H\"older, \eqref{small it}, and Sobolev embedding, on each $I_j\times\R$ we estimate
\begin{align*}
\|\langle \nabla\rangle u\|_{S^0(I_j)}
&\lesssim \|\langle \nabla\rangle u(t_j)\|_2 + \|\langle \nabla\rangle \bigl(|u|^{2p}u\bigr)\|_{L_t^{4/3}L_x^1}\\
&\lesssim \|u_0\|_{H^1(\R)} + \||u|^{2p}\|_{L_t^{2}L_x^{1}}\|\langle \nabla\rangle u\|_{L_t^4L_x^\infty}\\
&\lesssim \|u_0\|_{H^1(\R)} + \|u\|_{L_{t,x}^8}^\eps \|u\|_{L_t^{\frac{8(2p-\eps)}{4-\eps}}L_x^{\frac{8(2p-\eps)}{8-\eps}}}^{2p-\eps}\|\langle \nabla\rangle u\|_{S^0(I_j)}\\
&\lesssim \|u_0\|_{H^1(\R)} + \delta^\eps \bigl\||\nabla|^{\frac{8p-16-\eps}{8(2p-\eps)}}u\bigr\|_{S^0(I_j)}^{2p-\eps}\|\langle \nabla\rangle u\|_{S^0(I_j)}\\
&\lesssim \|u_0\|_{H^1(\R)} + \delta^\eps \|\langle \nabla\rangle u\|_{S^0(I_j)}^{2p+1 -\eps}.
\end{align*}
A standard continuity argument yields
$$
\|\langle \nabla\rangle u\|_{S^0(I_j)} \lesssim \|u_0\|_{H^1(\R)},
$$
provided $\delta$ is chosen sufficiently small depending on $\|u_0\|_{H^1(\R)}$.  Summing these bounds over all subintervals $I_j$ we derive
\begin{align}\label{sb}
\|\langle \nabla\rangle u\|_{S^0(I_j)} \leq C(\|u_0\|_{H^1(\R)}).
\end{align}

We now use \eqref{sb} to prove asymptotic completeness, that is, there exist unique $u_{\pm}$ such that
\begin{equation}\label{scat lim}
\|u(t)-e^{it\Delta}u_{\pm}\|_{H^1(\R)}\to 0 \quad \text{as }t\to\pm\infty.
\end{equation}
By time reversal symmetry, it suffices to prove the claim for positive times only.  For $t>0$, we define $v(t):=e^{-it\Delta}u(t)$.
We will  show that $v(t)$ converges in $H^1_x$ as $t\to +\infty$, and define $u_+$ to be the limit.

Indeed, by Duhamel's formula,
\begin{equation}\label{def v}
v(t)=u_0-i\int_0^t e^{-is\Delta} \bigl(|u|^{2p} u\bigr)(s)\,ds.
\end{equation}
Therefore, for $0<\tau<t$,
$$
v(t)-v(\tau)=-i\int_{\tau}^t e^{-is\Delta}\bigl(|u|^{2p}u\bigr)(s)\,ds.
$$
Arguing as above, by Lemma \ref{lemma linear strichartz} and Sobolev embedding,
\begin{align*}
\|v(t)-v(\tau)\|_{H^1(\R)}
&\lesssim \|\langle \nabla \rangle \bigl(|u|^{2p}u\bigr)\|_{L_t^{4/3}L_x^1([t,\tau]\times\R)}\\
&\lesssim \|u\|_{L_{t,x}^8([t,\tau]\times\R)}^\eps \|\langle \nabla \rangle u\|_{S^0([t,\tau])}^{2p+1-\eps}.
\end{align*}
Thus, by \eqref{mb} and \eqref{sb},
$$
\|v(t)-v(\tau)\|_{H^1(\R)}\to 0 \mbox{  as } \tau,t\to \infty.
$$
In particular, this implies $u_+$ is well defined and inspecting \eqref{def v} we find
$$
u_+=u_0-i\int_0^{\infty}e^{-is\Delta}(|u|^{2p}u)(s)\, ds.
$$
Using the same estimates as above, it is now an easy matter to derive \eqref{scat lim}.  This completes the proof of Theorem~\ref{scat}.

\section{Almost conservation law}\label{S:ac}
As mentioned in the introduction, in order to prove global well-posedness for \eqref{nls} it suffices to obtain \emph{a priori} control
over the  $H_x^s$ norm of solutions to \eqref{nls}.  However, the $H_x^s$ norm is not a conserved quantity.  Nevertheless, it can be
controlled by the $H_x^1$ norm of the modified solution $I_Nu$ (see \eqref{i3}).  While we do have conservation of energy for \eqref{nls},
$I_Nu$ is not a solution to \eqref{nls} and hence we expect an energy increment.  In this section, we prove that the energy increment is
small on intervals where the Morawetz norm is small, thus transfering the problem to controlling the Morawetz norm globally.

\begin{proposition}[Energy increment]\label{ac law}
Let $s>\frac{k-2}{2k-1}$ and let $u$ be an $H_x^s$ solution to \eqref{nls} on the spacetime slab $[t_0,T]\times \R$ with $E(I_N u(t_0))\le 1$.
Suppose in addition that
\begin{align}\label{ms}
\|u\|_{\ma([t_0,T]\times \R)}\le \eta
\end{align}
for a sufficiently small $\eta>0$ (depending on $k$ and on $E(I_N u(t_0))$). Then, for $N$ sufficiently large (depending on $k$ and on $E(I_N u(t_0))$),
\begin{equation}\label{eng growth}
\sup_{t\in [t_0,T]}E(I_N u(t))=E(I_N u(t_0))+ N^{-1+}.
\end{equation}
\end{proposition}

\begin{proof}
Fix $t\in[t_0,T]$ and define
$$
\|u\|_{Z(t)}: =\|\nabla P_{\leq 1}  u\|_{S^0([t_0,t])}
  + \sup_{(q,r) \ \text{admissible}} \Bigl(\sum_{N> 1}\|\nabla P_N u\|^2_{L_t^qL_x^r([t_0,t]\times\R)}\Bigr)^{1/2}.
$$

We observe the inequality
\begin{equation}\label{square sum}
\Bigl\|\Bigl(\sum_{N\in 2^\Z} |f_N|^2 \Bigr)^{1/2}\Bigr\|_{L_t^qL_x^r} \leq \Bigl(\sum_{N\in 2^\Z} \|f_N\|^2_{L_t^qL_x^r}\Bigr)^{1/2}
\end{equation}
for all $2 \leq q,r \leq \infty$ and arbitrary functions $f_N$, which one proves by interpolating between the trivial
cases $(2,2)$, $(2,\infty)$, $(\infty,2)$, and $(\infty,\infty)$. In particular, \eqref{square sum} holds for all
admissible exponents $(q,r)$.  Combining this with the Littlewood-Paley inequality, we find\footnote{Strictly speaking, as the Littlewood-Paley
square function is not bounded on $L_x^\infty$, the inequality does not hold for the Schr\"odinger-admissible pair $(4,\infty)$.  However,
this particular estimate will not be needed in the proof of Proposition~\ref{ac law} and we thus make the convention that in the proof of this
proposition alone the $S^0$ norm is the supremum over all admissible pairs except $(4,\infty)$.}
\begin{align*}
\| u \|_{L_t^qL_x^r}
\lesssim \Bigl\|\Bigl(\sum_{N\in 2^\Z} |P_N u|^2\Bigr)^{1/2}\Bigr\|_{L_t^qL_x^r}
\lesssim \Bigl(\sum_{N\in 2^\Z} \|P_N u \|^2_{L_t^qL_x^r}\Bigr)^{1/2}.
\end{align*}
In particular,
$$
\|\nabla u\|_{S^0([t_0,t])}\lesssim \|u\|_{Z(t)}.
$$
Moreover, using Lemma~\ref{lemma linear strichartz}, the fact that the Littlewood-Paley operators $P_N$ commute with $i\partial_t+\Delta$,
the Littlewood-Paley inequality, together with the dual of \eqref{square sum}, we get
\begin{align}\label{si}
\|u\|_{Z(t)}\lesssim \|u(t_0)\|_{\ho}+\|\nabla (iu_t+\Delta u)\|_{L^{q'}_tL^{r'}_{x}([t_0,t]\times\R)},
\end{align}
for any admissible pair $(q,r)$.

Now define
$$
\zit :=\|I_Nu\|_{Z(t)}.
$$

\begin{lemma}\label{zit lemma} Under the hypotheses of Proposition~\ref{ac law},
\begin{align}
Z_I(t)
&\lesssim \|\nabla I_Nu(t_0)\|_2+N^{-(k+2)}Z_I(t)^{2k+1}+\eta^{\frac {14k}{3k-1}}Z_I(t)^{1+\frac {2k(3k-8)}{3k-1}}\notag\\
&\quad +\eta^{\frac {16}3}\sup_{s\in [t_0,t]} E(I_Nu(s))^{\frac {3k-8}{3(k+1)}}Z_I(t). \label{zit control}
\end{align}
\end{lemma}

\begin{proof}
Throughout this proof, all spacetime norms are on $[t_0,t]\times \R$.  By \eqref{si} and H\"older's inequality, combined with the fact
that $\nabla I_N$ acts as a derivative (as the multiplier of $\nabla I_N$ is increasing in $|\xi|$), we estimate
\begin{align*}
Z_I(t)&\lesssim \|\nabla I_N u(t_0)\|_2+\|\nabla I_N (|u|^{2k}u)\|_{\frac 65,\frac 65}\\
&\lesssim \|\nabla I_N u(t_0)\|_2+\|u\|_{3k,3k}^{2k}\|\nabla I_N u\|_{6,6}\\
&\lesssim \|\nabla I_N u(t_0)\|_2+\|u\|_{3k,3k}^{2k}Z_I(t).
\end{align*}
To estimate $\|u\|_{3k,3k}$, we decompose $u:=u_{\le 1}+u_{1<\cdot \le N}+u_{>N}$.  To estimate the low frequencies, we use
interpolation, \eqref{ms}, Bernstein, and the fact that the operator $I_N$ is the identity on frequencies $|\xi|\leq 1$ to get
\begin{align*}
\|u_{\le 1}\|_{3k,3k}
&\lesssim \|u_{\le 1}\|_{8,8}^{\frac 8{3k}}\|u_{\le 1}\|_{\infty,\infty}^{1-\frac 8{3k}}\\
&\lesssim \eta^{\frac 8{3k}}\|u_{\le 1}\|_{\infty,2k+2}^{1-\frac 8{3k}}\\
&\lesssim \eta^{\frac 8{3k}}\sup_{s\in [t_0,t]} E(I_N u(s))^{\frac {3k-8}{3k(2k+2)}}.
\end{align*}
To estimate the medium frequencies, we use interpolation, \eqref{ms}, Sobolev embedding, Bernstein, and the fact that the operator $I_N$
is the identity on frequencies $|\xi|\leq N$
\begin{align*}
\|u_{1<\cdot \le N}\|_{3k,3k}
&\lesssim \|u_{1<\cdot\le N}\|_{8,8}^{\frac 7{3k-1}}\|u_{1<\cdot\le N}\|_{24k,24k}^{\frac {3k-8}{3k-1}}\\
&\lesssim \eta^{\frac 7{3k-1}}\||\nabla|^{\frac 12-\frac 1{8k}}u_{1<\cdot\le N}\|_{24k,\frac {12k}{6k-1}}^{\frac {3k-8}{3k-1}}\\
&\lesssim \eta^{\frac 7{3k-1}}Z_I(t)^{\frac {3k-8}{3k-1}}.
\end{align*}
To estimate the high frequencies, we use Sobolev embedding and Lemma~\ref{basic property}
\begin{align*}
\|u_{>N}\|_{3k,3k}
&\lesssim \||\nabla|^{\frac 12-\frac 1k}u_{>N}\|_{3k,\frac {6k}{3k-4}}\\
&\lesssim N^{-\frac 12-\frac 1k}\|\nabla I_N u_{>N}\|_{3k,\frac{6k}{3k-4}}\\
&\lesssim N^{-\frac 12-\frac 1k}Z_I(t).
\end{align*}
Putting everything together, we derive \eqref{zit control}.
\end{proof}

Next, we control the energy increment in terms of the size of the modified solution $I_Nu$.
\begin{lemma}\label{ei}
Under the hypotheses of Proposition~\ref{ac law},
\begin{align}
\bigl| & \sup_{s\in [t_0,t]} E(I_Nu(s))-E(I_Nu(t_0))\bigr|\label{energy increment} \\
&\lesssim N^{-1+}\Bigl(Z_I(t)^{2k+2}+ \eta^{\frac {16}{3}}Z_I(t)^2\sup_{s\in [t_0,t]} E(I_N u(s))^{\frac {3k-8}{3(k+1)}}\notag\\
&\qquad \qquad  + \sum_{J=3}^{2k+2}\eta^{\frac {4(2k+2-J)}{2k-1}}Z_I(t)^J \sup_{s\in [t_0,t]} E(I_Nu(s))^{\frac {(2k-5)(2k+2-J)}{(2k-1)(2k+2)}}\Bigr)\notag\\
&\quad +N^{-1+}\Bigl(\zit^{2k+1}+ \eta^{\frac {16}{3}}Z_I(t)\sup_{s\in [t_0,t]} E(I_N u(s))^{\frac {3k-8}{3(k+1)}}\Bigr)\notag\\
&\qquad \qquad\times \Bigl(\zit^{2k+1}+\eta^{\frac 43}\sup_{s\in [t_0,t]} E(I_N u(s))^{\frac {6k-1}{3(2k+2)}}\Bigr)\notag\\
&\quad +N^{-1+}\sum_{J=3}^{2k+2}\eta^{\frac{4(2k+2-J)}{2k-1}}\zit^{J-1}\sup_{s\in [t_0,t]} E(I_N u(s))^{\frac {(2k-5)(2k+2-J)}{(2k-1)(2k+2)}}\notag\\
&\qquad\qquad \times \Bigl(\zit^{2k+1}+\eta^{\frac 43}\sup_{s\in [t_0,t]} E(I_N u(s))^{\frac {6k-1}{3(2k+2)}}\Bigr).\notag
\end{align}
\end{lemma}

\begin{proof}
As
\begin{align*}
\frac d{dt}E(u(t))=\Re \int \bar u_t(|u|^{2k}u-\Delta u)\,dx=\Re \int \bar u_t(|u|^{2k}u-\Delta u-iu_t)\,dx,
\end{align*}
we obtain
\begin{align*}
\frac d{dt}E(Iu(t))&=\Re\int I\bar u_t (|Iu|^{2k}Iu-\Delta I u-iIu_t)\,dx\\
&=\Re \int I\bar u_t (|Iu|^{2k}Iu-I(|u|^{2k}u))\,dx.
\end{align*}
Using the Fundamental Theorem of Calculus and Plancherel, we write\footnote{Throughout this proof we use the abbreviation $m:=m_N$.}
\begin{align*}
E(Iu(t))&-E(Iu(t_0))\\
&=\Re\itxi \Bigl(\symb\Bigr)\\
&\qquad\qquad\qquad \widehat {\overline{I\partial_tu}}(\xi_1) \widehat{Iu}(\xi_2)\cdots\widehat{\overline{Iu}}(\xi_{2k+1}) \widehat{Iu}(\xi_{2k+2})\,d\sigma(\xi)\,ds.
\end{align*}
As $iu_t=-\Delta u+|u|^{2k}u$, we thus need to control
\begin{align}\label{term1}
\Bigl|\itxi&\Bigl(\symb\Bigr)\\
&\Delta \widehat{\overline {Iu}}(\xi_1)\widehat {Iu}(\xi_2)\cdots \widehat{\overline{Iu}}(\xi_{2k+1})\widehat {Iu}(\xi_{2k+2})\,d\sigma(\xi)\,ds\Bigr|\notag
\end{align}
and
\begin{align}\label{term2}
\Bigl|\itxi&\Bigl(\symb\Bigr)\\
&\widehat {\overline{I(|u|^{2k}u)}} (\xi_1)\widehat{Iu}(\xi_2)\cdots \widehat{\overline{Iu}}(\xi_{2k+1})\widehat {Iu}(\xi_{2k+2})\,d\sigma(\xi)\,ds\Bigr|.\notag
\end{align}

We first estimate \textbf{\eqref{term1}}.  To this end, we decompose
$$
u:=\sum_{N\geq 1}P_N u
$$
with the convention that $P_1 u: =P_{\leq 1}u$.  Using this notation and symmetry, we estimate
\begin{equation}
\eqref{term1}\lesssim \sum_{\substack{N_1, \dots, N_{2k+2}\geq 1 \\N_2\ge N_3\ge\cdots\geq N_{2k+2}}}B(N_1,\dots,N_{2k+2}),\label{al2}
\end{equation}
where
\begin{align*}
&B(N_1,\dots,N_{2k+2})\\
&\qquad:=\Bigl|\itxi\Bigl(\symb\Bigr)    \\
&\qquad \qquad \qquad \Delta\widehat{\overline{Iu_{N_1}}}(\xi_1)\widehat{Iu_{N_2}}(\xi_2)\cdots \widehat{\overline{Iu_{N_{2k+1}}}}(\xi_{2k+2}) \widehat {Iu_{N_{2k+2}}}(\xi_{2k+2})d\sigma(\xi)ds\Bigr|.
\end{align*}

\noindent{\bf Case $I$:} $N_1> 1$, $N_2\ge \cdots \ge N_{2k+2}> 1$.

\noindent \textbf{Case $I_a$:}  $N\gg N_2$.

In this case,
$$
m(\xi_2+\xi_3+\cdots+\xi_{2k+2})=m(\xi_2)=\dots=m(\xi_{2k+2})=1.
$$
Thus,
$$B(N_1,\dots,N_{2k+2})=0$$
and the contribution to the right-hand side of \eqref{al2} is zero.

\noindent \textbf{Case $I_b$:} $N_2\gtrsim N\gg N_3$.

As $\sum_{i=1}^{2k+2}\xi_i=0$, we must have $N_1\sim N_2$.  Thus, by the Fundamental Theorem of Calculus,
\begin{align*}
\Bigl|\symb\Bigr|&=\Bigl|1-\frac {m(\xi_2+\cdots+\xi_{2k+2})}{m(\xi_2)}\Bigr| \\
&\lesssim \Bigl|\frac {\nabla m(\xi_2)(\xi_3+\cdots+\xi_{2k+2})}{m(\xi_2)}\Bigr|\lesssim \frac {N_3}{N_2}.
\end{align*}
Applying the multilinear multiplier theorem (cf. \cite{cmfourier, cmaster}), Sobolev embedding, Bernstein, and recalling that $N_j>1$, we estimate
\begin{align*}
&B(N_1,\dots,N_{2k+2}) \\
&\qquad \lesssim \frac {N_3}{N_2}\|\Delta I u_{N_1}\|_{6,6}\|Iu_{N_2}\|_{6,6}\|Iu_{N_3}\|_{6,6}\prod_{j=4}^{2k+2}\|Iu_{N_j}\|_{2(2k-1),2(2k-1)}\\
&\qquad \lesssim \frac {N_1}{N_2^2}\prod_{j=1}^3\|\nabla I u_{N_j}\|_{6,6} \prod_{j=4}^{2k+2}\||\nabla|^{\frac {k-2}{2k-1}} Iu_{N_j}\|_{2(2k-1),\frac {2(2k-1)}{2k-3}}\\
&\qquad \lesssim \frac 1{N_2}Z_I(t)^{2k+2}
\lesssim N^{-1+}N_2^{0-}Z_I(t)^{2k+2}.
\end{align*}
The factor $N_2^{0-}$ allows us to sum in $N_1,N_2,\dots, N_{2k+2}$, this case contributing at most $N^{-1+}Z_I(t)^{2k+2}$ to the
right-hand side of \eqref{al2}.

\noindent \textbf{Case $I_c$:} $N_2\gg N_3\gtrsim N$.

As $\sum_{i=1}^{2k+2}\xi_i=0$, we must have $N_1\sim N_2$.  Thus, as $m$ is decreasing,
$$
\Bigl|\symb\Bigr|\lesssim \frac {m(\xi_1)}{m(\xi_2)\cdots m(\xi_{2k+2})}.
$$
Using again the multilinear multiplier theorem, Sobolev embedding, Bernstein, and the fact that $m(\xi)|\xi|^{\frac{k+1}{2k-1}}$ is increasing
for $s>\frac{k-2}{2k-1}$, we estimate
\begin{align*}
& \ B(N_1,\dots,N_{2k+2})   \\
&\lesssim \frac {m(N_1)}{m(N_2)\cdots m(N_{2k+2})}\frac {N_1}{N_2N_3} \prod_{j=1}^3\|\nabla I u_{N_j}\|_{6,6}\prod_{j=4}^{2k+2}\||\nabla|^{\frac {k-2}{2k-1}}Iu_{N_j}\|_ {2(2k-1),\frac {2(2k-1)}{2k-3}}\\
&\lesssim \frac 1{N_3m(N_3)\prod_{j=4}^{2k+2}m(N_j)N_j^{\frac {k+1}{2k-1}}}\prod_{j=1}^3\|\nabla I u_{N_j}\|_{6,6} \prod_{j=4}^{2k+2}\|\nabla I u_{N_j}\|_{2(2k-1),\frac{2(2k-1)}{2k-3}}\\
&\lesssim \frac 1{N_3m(N_3)}\|\nabla I u_{N_1}\|_{6,6} \|\nabla I u_{N_2}\|_{6,6}Z_I(t)^{2k}\\
&\lesssim N^{-1+}N_3^{0-}\|\nabla I u_{N_1}\|_{6,6} \|\nabla I u_{N_2}\|_{6,6}Z_I(t)^{2k}.
\end{align*}
The factor $N_3^{0-}$ allows us to sum over $N_3, \dots ,N_{2k+2}$.  To sum over $N_1$ and $N_2$, we use the fact that $N_1\sim N_2$
and Cauchy-Schwarz to estimate the contribution to the right-hand side of \eqref{al2} by
\begin{align*}
N^{-1+}\Bigl(\sum_{N_1>1}\|\nabla Iu_{N_1}\|_{6,6}^2\Bigr)^{\frac 12}\Bigl(\sum_{N_2>1}\|\nabla I u_{N_2}\|_{6,6}^2\Bigr)^{\frac 12}Z_I(t)^{2k}
\lesssim N^{-1+}Z_I(t)^{2k+2}.
\end{align*}

\noindent \textbf{Case $I_d$:} $N_2\sim N_3\gtrsim N$.

As $\sum_{i=1}^{2k+2}\xi_i=0$, we obtain $N_1\lesssim N_2$, and hence $m(N_1)\gtrsim m(N_2)$ and $m(N_1)N_1\lesssim m(N_2)N_2$.  Thus,
$$
\Bigl|\symb\Bigr|\lesssim \frac {m(N_1)}{m(N_2)m(N_3)\cdots m(N_{2k+2})}.
$$
Arguing as for Case $I_c$, we estimate
\begin{align*}
B(N_1,\dots, N_{2k+2})&\lesssim \frac {m(N_1)N_1}{m(N_2)N_2m(N_3)N_3\prod_{j=4}^{2k+2}m(N_j)N_j^{\frac
{k+1}{2k-1}}}Z_I(t)^{2k+2}\\
&\lesssim \frac 1{m(N_3)N_3}Z_I(t)^{2k+2}\\
&\lesssim N^{-1+}N_3^{0-}Z_I(t)^{2k+2}.
\end{align*}
The factor $N_3^{0-}$ allows us to sum over $N_1, \dots, N_{2k+2}$.  This case contributes at most $N^{-1+}Z_I(t)^{2k+2}$
to the right-hand side of \eqref{al2}.

\noindent \textbf{Case $II$:} There exists $1\leq j_0\leq 2k+2$ such that $N_{j_0}=1$.  Recall that by our convention, $P_1:=P_{\leq 1}$.

\noindent \textbf{Case $II_a$:} $N_1=1$.

Let $J$ be such that $N_2\ge \dots\geq N_J> 1=N_{J+1}=\dots=N_{2k+2}$.  Note that we may assume $J\geq 3$ since otherwise
$$
B(N_1,\dots,N_{2k+2})=0.
$$
Also, arguing as for Case $I_a$, if $N\gg N_2$ then
$$
B(N_1,\dots,N_{2k+2})=0.
$$
Thus, we may assume $N_2\gtrsim N$.  In this case we cannot have $N_2\gg N_3$ since it would contradict $\sum_{i=1}^{2k+2}\xi_i=0$ and $N_1=1$.
Hence, we must have
$$
N_2\sim N_3\gtrsim N.
$$
As
$$
\Bigl|\symb\Bigr|\lesssim \frac 1{m(N_2)m(N_3)\cdots m(N_{2k+2})},
$$
we use the multilinear multiplier theorem and Sobolev embedding to estimate
\begin{align*}
& \ B(N_1,\dots,N_{2k+2}) \\
&\lesssim \frac {N_1}{m(N_2)N_2m(N_3)N_3m(N_4)\cdots m(N_{2k+2})}\prod_{j=1}^3\|\nabla Iu_{N_j}\|_{6,6}\\
&\qquad \qquad \times \prod_{j=4}^J\||\nabla|^{\frac {k-2}{2k-1}}Iu_{N_j}\|_{2(2k-1),\frac {2(2k-1)}{2k-3}}\prod _{j=J+1}^{2k+2}\|Iu_{N_j}\|_{2(2k-1),2(2k-1)} \\
&\lesssim \frac 1{m(N_2)N_2m(N_3)N_3\prod_{j=4}^Jm(N_j)N_j^{\frac{k+1}{2k-1}}}Z_I(t)^J \prod_{j=J+1}^{2k+2}\|I u_{N_j}\|_{2(2k-1),2(2k-1)}\\
&\lesssim N^{-2+}N_2^{0-}Z_I(t)^J\prod_{j=J+1}^{2k+2}\|Iu_{N_j}\|_{2(2k-1),2(2k-1)}.
\end{align*}
Applying interpolation, \eqref{ms}, and Bernstein, we bound
\begin{align}
\|Iu_{\leq 1}\|_{2(2k-1),2(2k-1)}
&\lesssim \|Iu_{\leq 1}\|_{8,8}^{\frac4{2k-1}}\| I u_{\leq 1}\|_{\infty,\infty}^{\frac {2k-5}{2k-1}}\label{low}\\
&\lesssim \eta^{\frac 4{2k-1}}\sup_{s\in[t_0,t]} E(Iu(s))^{\frac {2k-5}{(2k-1)(2k+2)}}. \nonumber
\end{align}
Thus,
\begin{align*}
B(N_1, \dots,  N_{2k+2})
\lesssim N^{-2+} N_2^{0-}\eta^{\frac {4(2k+2-J)}{2k-1}} Z_I(t)^J \sup_{s\in[t_0,t]}E(I u(s))^{\frac {(2k-5)(2k+2-J)}{(2k-1)(2k+2)}}.
\end{align*}
The factor $N_2^{0-}$ allows us to sum in $N_2, \dots, N_J$.  This case contributes at most
$$
N^{-2+} \sum_{J=3}^{2k+2}\eta^{\frac {4(2k+2-J)}{2k-1}}Z_I(t)^J \sup_{s\in[t_0,t]}E(I u(s))^{\frac {(2k-5)(2k+2-J)}{(2k-1)(2k+2)}}
$$
to the right-hand side of \eqref{al2}.

\noindent \textbf{Case $II_b$:} $N_1> 1$ and $N_2=\dots=N_{2k+2}=1$.

As $\sum_{i=1}^{2k+2}\xi_i=0$, we obtain $N_1\lesssim 1$ and thus, taking $N$ sufficiently large depending on $k$, we get
$$
\symb=0.
$$
This case contributes zero to the right-hand side of \eqref{al2}.

\noindent \textbf{Case $II_c$:} $N_1> 1$ and $N_2>1=N_3=\dots=N_{2k+2}$.

As $\sum_{i=1}^{2k+2}\xi_i=0$, we must have $N_1\sim N_2$.  If $N_1\sim N_2\ll N$, then
$$
\symb=0
$$
and the contribution is zero.  Thus, we may assume $N_1\sim N_2 \gtrsim N$.

Applying the Fundamental Theorem of Calculus,
\begin{align*}
\Bigl|\symb\Bigr|&=\Bigl|1-\frac {m(\xi_2+\cdots+\xi_{2k+2})}{m(\xi_2)}\Bigr| \\
&\lesssim \Bigl|\frac {\nabla m(\xi_2)}{m(\xi_2)}\Bigr|\lesssim \frac {1}{N_2}.
\end{align*}
By the multilinear multiplier theorem,
\begin{align*}
B(N_1,\dots,N_{2k+2})
&\lesssim \frac {1}{N_2}\|\Delta I u_{N_1}\|_{6,6}\|Iu_{N_2}\|_{6,6}\prod_{j=3}^{2k+2}\|Iu_{N_j}\|_{3k,3k}\\
&\lesssim \frac {N_1}{N_2^2}\|\nabla I u_{N_1}\|_{6,6} \|\nabla I u_{N_2}\|_{6,6} \|Iu_{\leq 1}\|_{3k,3k}^{2k}\\
&\lesssim N^{-1+}N_2^{0-}Z_I(t)^2\|Iu_{\leq 1}\|_{3k,3k}^{2k}.
\end{align*}
The factor $N_2^{0-}$ allows us to sum in $N_1$ and $N_2$.  Using interpolation, \eqref{i1}, \eqref{ms}, and Bernstein, we estimate
\begin{align*}
\|Iu_{\le 1}\|_{3k,3k}
&\lesssim \|Iu_{\le 1}\|_{8,8}^{\frac 8{3k}}\|Iu_{\le 1}\|_{\infty,\infty}^{1-\frac 8{3k}}\\
&\lesssim \eta^{\frac 8{3k}}\|Iu_{\le 1}\|_{\infty,2k+2}^{1-\frac 8{3k}}\\
&\lesssim \eta^{\frac 8{3k}}\sup_{s\in [t_0,t]} E(I u(s))^{\frac {3k-8}{3k(2k+2)}}.
\end{align*}
Thus, this case contributes at most
$$
N^{-1+}\eta^{\frac {16}{3}}Z_I(t)^2\sup_{s\in [t_0,t]} E(I u(s))^{\frac {3k-8}{3(k+1)}}
$$
to the right-hand side of \eqref{al2}.

\noindent \textbf{Case $II_d$:} $N_1> 1$ and there exists $J\ge 3$ such that $N_2\ge \dots \ge N_J> 1 =N_{J+1}=\dots=N_{2k+2}$.

To estimate the contribution of this case, we argue as for Case $I$; the only new ingredient is that the low frequencies are estimated via \eqref{low}.
This case contributes at most
\begin{align*}
N^{-1+}\sum_{J=3}^{2k+2}\eta^{\frac {4(2k+2-J)}{2k-1}}Z_I(t)^J \sup_{s\in[t_0,t]} E(Iu(s))^{\frac {(2k-5)(2k+2-J)}{(2k-1)(2k+2)}}
\end{align*}
to the right-hand side of \eqref{al2}.

Putting everything together, we get
\begin{align}
\eqref{term1}
&\lesssim N^{-1+}Z_I(t)^{2k+2} + N^{-1+}\eta^{\frac {16}{3}}Z_I(t)^2\sup_{s\in [t_0,t]} E(I u(s))^{\frac {3k-8}{3(k+1)}} \notag\\
&\quad + N^{-1+}\sum_{J=3}^{2k+2}\eta^{\frac {4(2k+2-J)}{2k-1}}Z_I(t)^J \sup_{s\in [t_0,t]} E(Iu(s))^{\frac {(2k-5)(2k+2-J)}{(2k-1)(2k+2)}}.\label{term1 est}
\end{align}

We turn now to estimating \textbf{\eqref{term2}}.  Again we decompose
$$
u:=\sum_{N\geq 1}P_N u
$$
with the convention that $P_1 u: =P_{\leq 1}u$.  Using this notation and symmetry, we estimate
\begin{align*}
\eqref{term2}\lesssim \sum_{\substack{N_1, \dots, N_{2k+2}\geq 1 \\ N_2\ge \cdots\ge N_{2k+2}}} C(N_1,\cdots,N_{2k+2}),
\end{align*}
where
\begin{align*}
C(& N_1,\cdots, N_{2k+2})\\
&:=\Bigl|\itxi\Bigl(\symb\Bigr)\\
&\qquad \widehat{\overline{P_{N_1}I(|u|^{2k}u)}}(\xi_1)\widehat{I u_{N_2}}(\xi_2)\cdots \widehat{\overline{Iu_{N_{2k+1}}}}
(\xi_{2k+1})\widehat{Iu_{N_{2k+2}}}(\xi_{2k+2})\, d\sigma(\xi)\,ds\Bigr|.
\end{align*}
In order to estimate $C(N_1,\cdots, N_{2k+2})$ we make the observation that in estimating $B(N_1,\cdots, N_{2k+2})$,
for the term involving the $N_1$ frequency we only used the bound
\begin{align}\label{1use}
\|P_{N_1}I\Delta u\|_{6,6}\lesssim N_1 \|\nabla Iu_{N_1}\|_{6,6}\lesssim N_1 Z_I(t).
\end{align}
Thus, to estimate \eqref{term2} it suffices to prove
\begin{equation}\label{al3}
\|P_{N_1}I(|u|^{2k}u)\|_{6,6}\lesssim Z_I(t)^{2k+1}+\eta^{\frac 43}\sup_{s\in [t_0,t]} E(Iu(s))^{\frac {6k-1}{3(2k+2)}},
\end{equation}
for then, arguing as for \eqref{term1} and substituting \eqref{al3} for \eqref{1use}, we obtain
\begin{align*}
\eqref{term2}
&\lesssim N^{-1+}\Bigl(\zit^{2k+1}+ \eta^{\frac {16}{3}}Z_I(t)\sup_{s\in [t_0,t]} E(I u(s))^{\frac {3k-8}{3(k+1)}}\Bigr)\\
&\qquad \qquad\times \Bigl(\zit^{2k+1}+\eta^{\frac 43}\sup_{s\in [t_0,t]} E(I u(s))^{\frac {6k-1}{3(2k+2)}}\Bigr)\\
&\quad +N^{-1+}\sum_{J=3}^{2k+2}\eta^{\frac{4(2k+2-J)}{2k-1}}\zit^{J-1}\sup_{s\in [t_0,t]} E(I u(s))^{\frac {(2k-5)(2k+2-J)}{(2k-1)(2k+2)}}\\
&\qquad\qquad \times \Bigl(\zit^{2k+1}+\eta^{\frac 43}\sup_{s\in [t_0,t]} E(I u(s))^{\frac {6k-1}{3(2k+2)}}\Bigr).
\end{align*}
Thus, we are left to proving \eqref{al3}.  Using \eqref{i1} and the boundedness of the Littlewood-Paley operators, and decomposing
$u:=u_{\leq 1}+u_{>1}$, we estimate
\begin{align*}
\|P_{N_1}I(|u|^{2k}u)\|_{6,6}
&\lesssim \|u\|_{6(2k+1),6(2k+1)}^{2k+1}\\
&\lesssim\|u_{\le 1}\|_{6(2k+1),6(2k+1)}^{2k+1}+\|u_{>1}\|_{6(2k+1),6(2k+1)}^{2k+1}.
\end{align*}
Applying interpolation, \eqref{ms}, and Bernstein, we estimate
\begin{align*}
\|u_{\le 1}\|_{6(2k+1),6(2k+1)}
&\lesssim\|u_{\le 1}\|_{8,8}^{\frac 4{3(2k+1)}}\|u_{\le 1}\|_{\infty,\infty}^{\frac{6k-1}{3(2k+1)}}\\
&\lesssim \eta^ {\frac4{3(2k+1)}}\sup_{s\in [t_0,t]} E(I u(s))^{\frac{6k-1}{3(2k+2)(2k+1)}}.
\end{align*}
Finally, by Sobolev embedding and \eqref{i2},
\begin{align*}
\|u_{>1}\|_{6(2k+1),6(2k+1)}
&\lesssim\||\nabla|^{\frac12-\frac1{2(2k+1)}}u_{>1}\|_{6(2k+1),\frac{6(2k+1)}{6k+1}}
\lesssim \zit.
\end{align*}
Putting things together, we derive \eqref{al3}.

This completes the proof of Lemma~\ref{ei}.
\end{proof}

Next, we combine Lemmas~\ref{zit lemma} and \ref{ei} to derive Proposition~\ref{ac law}.  Indeed, Proposition~\ref{ac law}
follows immediately from Lemmas~\ref{zit lemma} and \ref{ei}, if we establish
\begin{equation*}
Z_I(t)\lesssim 1 \quad \text{and} \quad \sup_{s\in [t_0,t]} E(I_N u(s))\lesssim 1 \quad \text{for all } t\in[t_0,T].
\end{equation*}
As by assumption $E(I_N u(t_0))\leq 1$, it suffices to show that
\begin{equation}\label{bdd1}
Z_I(t)\lesssim\|\nabla I_N u(t_0)\|_2 \quad \text{for all } t\in[t_0,T]
\end{equation}
and
\begin{equation}\label{bdd2}
\sup_{s\in [t_0,t]} E(I_N u(s))\lesssim E(I_N u(t_0)) \quad \text{for all } t\in[t_0,T].
\end{equation}
We achieve this via a bootstrap argument. Let
\begin{align*}
\Omega_1&:=\{t\in[t_0,T]:\, Z_I(t)\le C_1\|\nabla I_N u(t_0)\|_2,\\
    & \qquad \qquad \qquad \quad \sup_{s\in [t_0,t]} E(I_N u(s))\le C_2 E(I_N u(t_0))\}\\
\Omega_2&:=\{t\in[t_0,T]:\, Z_I(t)\le 2 C_1\|\nabla I_N u(t_0)\|_2,\\
    &\qquad \qquad \qquad  \quad \sup_{s\in [t_0,t]} E(I_N u(s)) \le 2C_2 E(I_Nu(t_0))\}.
\end{align*}
In order to run the bootstrap argument successfully, we need to check four things:

$\bullet$ $\Omega_1\neq \emptyset$.  This is satisfied as $t_0\in \Omega_1$ if we take $C_1$ and $C_2$ sufficiently large.

$\bullet$ $\Omega_1$ is a closed set.  This follows from Fatou's Lemma.

$\bullet$ If $t\in \Omega_1$, then there exists $\eps>0$ such that $[t,t+\eps]\in \Omega_2$.  This follows from the Dominated Convergence Theorem
combined with \eqref{zit control} and \eqref{energy increment}.

$\bullet$ $\Omega_2\subset \Omega_1$.  This follows from \eqref{zit control} and \eqref{energy increment} taking $C_1$ and $C_2$
sufficiently large depending on absolute constants (like the Strichartz constant) and choosing $N$ sufficiently
large and $\eta$ sufficiently small depending on $C_1$, $C_2$, $k$, and $E(I_N u(t_0))$.

This finally proves Proposition \ref{ac law}.

\section{Proof of Theorem \ref{main}}\label{S:boot}
Given Proposition \ref{ac law}, the proof of global well-posedness for \eqref{nls} is reduced to showing
\begin{equation}\label{ma bound}
\|u\|_{\ma(\R\times\R)}\le C(\|u_0\|_{\hs}).
\end{equation}
This also implies scattering, as we will see later.

By Proposition~\ref{fpim},
\begin{align}\label{ma control}
\|u\|_{\ma(\ir)}\lesssim \|u_0\|_2^{3/4} \|u\|_{L_t^\infty \dot H^{1/2}_x (\ir)}^{1/4}
\end{align}
on any spacetime slab $\ir$ on which the solution to \eqref{nls} exists and lies in $H_x^{1/2}$.
However, the $H_x^{1/2}$ norm of the solution is not a conserved quantity either, and in order to control it we must resort to
the $\hs$ bound on the solution.  Thus, in order to obtain a global Morawetz estimate, we need a global $\hs$ bound.
This sets us up for a bootstrap argument.

Let $u$ be the solution to \eqref{nls}. As $E(I_N u_0)$ is not necessarily small, we first rescale the solution such that the energy
of the rescaled initial data satisfies the conditions in Proposition~\ref{ac law}. By scaling,
$$
\ulam(x,t):=\lambda^{-\frac 1k} u({\lambda}^{-2} t,{\lambda^{-1}} x)
$$
is also a solution to \eqref{nls} with initial data
$$
u_0^{\lambda}(x):=\lambda^{-\frac1k}u_0({\lambda}^{-1} x).
$$
By \eqref{i4} and Sobolev embedding,
\begin{align*}
\|\nabla I_N u_0^{\lambda}\|_2
&\lesssim N^{1-s}\|u_0^{\lambda}\|_{\dhs}=N^{1-s}\lambda^{\frac 12-\frac 1k-s}\|u_0\|_{\dhs},\\
\|I_N u_0^{\lambda}\|_{2k+2}&\lesssim \|u_0^{\lambda}\|_{2k+2}=\lambda^{\frac 1{2k+2}-\frac 1k}\|u_0\|_{2k+2}\lesssim\lambda^{\frac 1{2k+2}-\frac 1k}\|u_0\|_{\hs}.
\end{align*}
As $s>\frac 12-\frac 1k$, choosing $\lambda$ sufficiently large (depending on $\|u_0\|_{\hs}$ and $N$) such that
\begin{equation}\label{lambdachoice}
N^{1-s} \lambda^{\frac{1}{2} -\frac{1}{k} -s} \|u_0\|_{\hs} \ll 1 \quad \text{and}\quad \lambda^{\frac 1{2k+2}-\frac 1k}\|u_0\|_{\hs}\ll 1,
\end{equation}
we get
$$
E(I_N u_0^{\lambda})\ll1.
$$
We now show that there exists an absolute constant $C_1$ such that
\begin{equation}\label{rescaled ma}
\|\ulam\|_{\ma(\R\times\R)}\le C_1\lambda^{\frac 78(\frac12-\frac 1k)}.
\end{equation}
Undoing the scaling, this yields \eqref{ma bound}.

We prove \eqref{rescaled ma} via a bootstrap argument.  By time reversal symmetry, it suffices to argue for positive times only.  Define
\begin{align*}
\Omega_1&:=\{t\in[0,\infty):\, \|\ulam\|_{\ma([0,t]\times\R)}\le C_1\lambda^{\frac 78(\frac12-\frac 1k)}\}, \\
\Omega_2&:=\{t\in[0,\infty):\, \|\ulam\|_{\ma([0,t]\times\R)}\le 2C_1\lambda^{\frac 78(\frac12-\frac 1k)}\}.
\end{align*}
In order to run the bootstrap argument, we need to verify four things:

1) $\Omega_1\neq\emptyset$.  This is obvious as $0\in\Omega_1$.

2) $\Omega_1$ is closed.  This follows from Fatou's Lemma.

3) $\Omega_2\subset\Omega_1$.

4) If $T\in\Omega_1$, then there exists $\eps>0$ such that $[T,T+\eps)\subset \Omega_2$.  This is a consequence of the local
well-posedness theory and the proof of 3). We skip the details.

Thus, we need to prove 3).  Fix $T\in \Omega_2$; we will show that in fact, $T\in\Omega_1$.
By \eqref{ma control} and the conservation of mass,
\begin{align*}
\|\ulam\|_{\ma([0,T]\times \R)}
&\lesssim \|u_0^{\lambda}\|_2^{\frac 34}\|\ulam\|^{\frac 14}_{L_t^\infty \dot H_x^{1/2}([0,T]\times\R)}\\
&\lesssim \lambda^{\frac 34(\tk)} C(\|u_0\|_2) \|\ulam\|_{L_t^\infty \dot H_x^{1/2}([0,T]\times\R)}^{\frac 14}.
\end{align*}
To control the factor $\|\ulam\|_{L_t^\infty \dot H_x^{1/2}([0,T]\times\R)}$, we decompose
$$
\ulam(t):=P_{\le N}\ulam(t)+P_{>N}\ulam(t).
$$
To estimate the low frequencies, we interpolate between the $\lxt$ norm and the $\ho$ norm and use the fact that $I_N$ is
the identity on frequencies $|\xi|\le N$
\begin{align*}
\|P_{\le N}\ulam(t)\|_{\dot H_x^{1/2}}
&\lesssim \|P_{\le N}\ulam(t)\|_2^{\frac 12}\|P_{\le N}\ulam(t)\|_{\ho}^{\frac 12}\\
&\lesssim\lambda^{\frac 12(\tk)}C(\|u_0\|_2)\|I_N \ulam(t)\|_{\ho}^{\frac 12}.
\end{align*}
To control the high frequencies, we interpolate between the $\lxt$ norm and the $\dhs$ norm and use Lemma \ref{basic property}
\begin{align*}
\|P_{>N}\ulam(t)\|_{\dot H^{1/2}_x}
&\lesssim \|P_{>N}\ulam(t)\|_{\lxt}^{1-\frac 1{2s}}\|P_{>N}\ulam(t)\|_{\dhs}^{\frac 1{2s}}\\
&\lesssim \lambda^{(1-\frac 1{2s})(\tk)}N^{\frac{s-1}{2s}}\|I_N\ulam(t)\|_{\ho}^{\frac1{2s}}\\
&\lesssim \lambda^{\frac 12(\tk)}\|I_N\ulam(t)\|_{\ho}^{\frac 1{2s}}.
\end{align*}
Collecting all these estimates, we get
\begin{align*}
\|\ulam\|_{\ma(\ztr)}
&\lesssim \lambda^{\frac 78(\frac12-\frac 1k)} C(\|u_0\|_2)
   \sup_{t\in [0,T]}\bigl(\|\nabla I_N \ulam(t)\|_2^{\frac 18}+\|\nabla I_N\ulam(t)\|_2^{\frac 1{8s}}\bigr).
\end{align*}
Thus, taking $C_1$ sufficiently large depending on $\|u_0\|_2$, we obtain $T\in \Omega_1$, provided
\begin{equation}\label{bdd kinetic}
\sup_{t\in [0,T]}\|\nabla I_N\ulam(t)\|_2\le 1.
\end{equation}
We now prove that $T\in\Omega_2$ implies \eqref{bdd kinetic}.  Indeed, let $\eta>0$ be a sufficiently small constant like in
Proposition~\ref{ac law} and divide $[0,T]$ into
$$
L\sim \biggl(\frac{\lambda^{\frac 78 (\frac 12 -\frac 1k)}}{\eta}\biggr)^8
$$
subintervals $I_j=[t_j,t_{j+1}]$ such that,
$$
\|\ulam\|_{\ma(I_j\times\R)}\le \eta.
$$
Applying Proposition \ref{ac law} on each of the subintervals $I_j$,
we get
$$
\sup_{t\in [0,T]} E(I_N\ulam(t))\le E(I_Nu_0^{\lambda})+ E(I_N u_0^{\lambda})LN^{-1+}.
$$
To maintain small energy during the iteration, we need
\begin{align*}
LN^{-1+}\sim\lambda^{7(\frac 12-\frac 1k)}N^{-1+}\ll 1,
\end{align*}
which combined with \eqref{lambdachoice} leads to
$$
\biggl(N^{\frac{1-s}{s+\frac 1k-\frac 12}}\biggr)^{7(\tk)}N^{-1+}\le c(\|u_0\|_{\hs})\ll 1.
$$
This may be ensured by taking $N$ large enough (depending only on $k$ and $\|u_0\|_{H^s(\R^)}$), provided that
$$
s>s(k):=\frac{8k-16}{9k-14}.
$$
As can be easily seen, $s(k)\to \frac 89$ as $k\to \infty$.

This completes the bootstrap argument and hence \eqref{rescaled ma}, and moreover \eqref{ma bound}, follow.  Therefore \eqref{bdd kinetic}
holds for all $T\in\R$ and the conservation of mass and Lemma~\ref{basic property} imply
\begin{align*}
\|u(T)\|_{\hs}&\lesssim \|u_0\|_{\lxt}+\|u(T)\|_{\dhs}\\
&\lesssim\|u_0\|_{\lxt}+\lambda^{s-(\tk)}\|\ulam(\lambda^2T)\|_{\dhs}\\
&\lesssim \|u_0\|_{\lxt}+\lambda^{s-(\tk)}\|I_N\ulam(\lambda^2T)\|_{H_x^1}\\
&\lesssim\|u_0\|_{\lxt}+\lambda^{s-(\tk)}(\|\ulam(\lambda^2T)\|_{\lxt}+\|\nabla I_N\ulam(\lambda^2T)\|_{\lxt})\\
&\lesssim\|u_0\|_{\lxt}+\lambda^{s-(\tk)}(\lambda^{\tk}\|u_0\|_{\lxt}+1)\\
&\lesssim  C(\|u_0\|_{\hs})
\end{align*}
for all $T\in \R$.  Hence,
\begin{align}\label{hsbdd}
\|u\|_{L_t^{\infty}\hs}\le C(\|u_0\|_{\hs}).
\end{align}

Finally, we prove that scattering holds in $\hs$ for $s>s_k$.  As the construction of the wave operators is standard (see \cite{cazbook}),
we content ourselves with proving asymptotic completeness.

The first step is to upgrade the global Morawetz estimate to global Strichartz control.
Let $u$ be a global $\hs$ solution to \eqref{nls}.  Then $u$ satisfies \eqref{ma bound}.
Let $\delta>0$ be a small constant to be chosen momentarily and split $\R$ into $L=L(\|u_0\|_{\hs})$ subintervals $I_j=[t_j, t_{j+1}]$ such that
$$
\|u\|_{\ma(I_j\times\R)}\le \delta.
$$
By Lemma~\ref{lemma linear strichartz}, \eqref{hsbdd}, and the fractional chain rule, \cite{christweinstein}, we estimate
\begin{align*}
\|\langle\nabla\rangle^s u\|_{S^0(I_j)}
&\lesssim \|u(t_j)\|_{\hs}+\|\langle\nabla\rangle^s\bigl(|u|^{2k}u\bigr)\|_{L_{t,x}^{6/5}(I_j\times\R)}\\
&\lesssim C(\|u_0\|_{\hs})+\|u\|_{L_{t,x}^{3k}}^{2k}\|\langle \nabla\rangle^s u\|_{L_{t,x}^6(I_j\times\R)},
\end{align*}
while by H\"older and Sobolev embedding,
\begin{align*}
\|u\|_{L_{t,x}^{3k}(I_j\times\R)}
&\lesssim \|u\|_{\ma(I_j\times\R)}^{\frac 7{3k-1}}\|u\|_{L_{t,x}^{24k}(I_j\times\R)}^{\frac{3k-8}{3k-1}}\\
&\lesssim \delta^{\frac 7{3k-1}}\||\nabla|^{\frac 12-\frac 1{8k}}u\|_{L_t^{24k}L_x^{\frac {12k}{6k-1}}(I_j\times\R)}\\
&\lesssim \delta^{\frac 7{3k-1}}\|\langle\nabla\rangle^s u\|_{S^0(I_j)}^{\frac{3k-8}{3k-1}}.
\end{align*}
Therefore,
$$
\|\langle\nabla\rangle^s u\|_{S^0(I_j)}
\lesssim C(\|u_0\|_{\hs})+\delta^{\frac {14k}{3k-1}}\|\langle\nabla\rangle^s u\|_{S^0(I_j)}^{1+\frac {2k(3k-8)}{3k-1}}.
$$
A standard continuity argument yields
$$
\|\langle\nabla\rangle^su\|_{S^0(I_j)}\le C(\|u_0\|_{\hs}),
$$
provided we choose $\delta$ sufficiently small depending on $k$ and $\|u_0\|_{\hs}$.  Summing over all subintervals $I_j$, we obtain
\begin{equation}\label{s bound}
\|\langle\nabla\rangle^s u\|_{S^0(\R)}\le C(\|u_0\|_{\hs}).
\end{equation}

We now use \eqref{s bound} to prove asymptotic completeness, that is, there exist unique $u_{\pm}$ such that
\begin{equation}\label{limit}
\lim_{t\to \pm\infty}\|u(t)-e^{it\Delta}u_{\pm}\|_{\hs}=0.
\end{equation}
Arguing as in Section~\ref{S:scat}, it suffices to see that
\begin{align}\label{lg}
\Bigl\|\int_{t}^\infty e^{-is\Delta}\bigl(|u|^{2k}u\bigr)(s)\,ds\Bigr\|_{\hs}\to 0 \quad \text{as } t\to \infty.
\end{align}
The estimates above yield
\begin{align*}
\Bigl\|\int_{t}^\infty e^{-is\Delta}\bigl(|u|^{2k}u\bigr)(s)\,ds\Bigr\|_{\hs}
\lesssim \|u\|_{\ma([t,\infty]\times\R)}^{\frac {14k}{3k-1}}\|\langle \nabla \rangle^s u\|_{S^0([t,\infty]\times\R)}^{1+\frac {2k(3k-8)}{3k-1}}.
\end{align*}
Using \eqref{ma bound} and \eqref{s bound} we derive \eqref{lg}.

This concludes the proof of Theorem~\ref{main}.

\end{proof}

\end{document}